\documentclass[12pt,a4paper]{article}%
\usepackage{amssymb}
\usepackage{amsmath}
\usepackage{amsfonts}
\usepackage{graphicx}%
\providecommand{\U}[1]{\protect\rule{.1in}{.1in}}

\begin{document}

\begin{center}
{\LARGE On Zermelo's theorem}\footnote[1]{The idea of this note emerged when
the authors were participating in the Trimester Program "Stochastic Dynamics
in Economics and Finance" at the Hausdorff Research Institute for Mathematics
(Bonn, May -- August 2013). Our analysis was inspired by Yurii Khomskii's
lectures on the methods of mathematical logic in game theory that were
delivered in the framework of the Trimester. We are grateful to Wolfgang
L\"{u}ck, Director of the Hausdorff Institute, for his kind invitation to the
Trimester Program. We thank Yurii Khomskii for fruitful discussions and
helpful comments.}\bigskip

Rabah Amir\footnote[2]{Department of Economics, University of Iowa, Iowa City,
IA 52242-1994, USA, rabah-amir@uiowa.edu .} and Igor V. Evstigneev\bigskip
\footnote[3]{Department of Economics, University of Manchester, Oxford Road,
M13 9PL, UK, igor.evstigneev@manchester.ac.uk (corresponding author).}\bigskip
\end{center}

\begin{quotation}
\textbf{Abstract.}{\small \ A famous result in game theory known as Zermelo's
theorem says that "in chess either White can force a win, or Black can force a
win, or both sides can force at least a draw". The present paper extends this
result to the class of all finite-stage two-player games of complete
information with alternating moves. It is shown that in any such game either
the first player has a winning strategy, or the second player has a winning
strategy, or both have unbeatable strategies.\bigskip}


\textbf{Key words:}{\small \ games of complete information, combinatorial game
theory, dynamic games, determinacy, Zermelo's theorem\bigskip}
\end{quotation}

\textbf{1.} In this note we generalize the following proposition usually
referred to in the modern game-theoretic literature as Zermelo's theorem (see
the paper [6] by Zermelo and its discussion and further references in [5]).

\textbf{Theorem 1.} \textit{In chess either White has a winning strategy, or
Black has a winning strategy, or both have strategies guaranteeing at least a
draw.}

The main result of this paper is as follows.

\textbf{Theorem 2.}\textit{\ In any finite-stage two-player game with
alternating moves, either (i) player 1 has a winning strategy, or (ii) player
2 has a winning strategy, or (iii) both have unbeatable strategies.}

We emphasize that in this theorem we speak of \textit{any }finite-stage
two-player game with alternating moves, with \textit{any }action sets and
\textit{any }real-valued payoff functions, so that the formulation of the
result has a maximum level of generality.

The notions of unbeatable and winning strategies are defined in the general
context as follows (cf. [1,4]). Consider a game of two players who select
strategies $\xi$ and $\eta$ from some sets and get payoffs $U(\xi,\eta)$ and
$V(\xi,\eta)$. We call a strategy $\xi$ of player 1 \textit{unbeatable (}resp.
\textit{winning)} if $U(\xi,\eta)\geq V(\xi,\eta)$ (resp. $U(\xi,\eta
)>V(\xi,\eta)$) for any $\eta$. Unbeatable and winning strategies $\eta$ of
player 2 are defined analogously in terms of the inequalities $U(\xi,\eta)\leq
V(\xi,\eta)$ and $U(\xi,\eta)<V(\xi,\eta)$ holding for each $\xi$.

The notions introduced are regarded in this work as primitive game solution
concepts. We do not try to reduce them to the conventional ones: saddle point,
Nash equilibrium, dominant strategy or their versions. This approach,
contrasting with the common game-theoretic methodology, constitutes the main
element of novelty in this paper. It enables us to obtain the above result in
its most general and natural form.

\textbf{2. }In a \textit{finite-stage game with alternating moves}, players 1
and 2 make moves (take actions) sequentially by selecting elements from two
given \textit{action sets} $A$ and $B$. At stage $0$ player 1 makes a move
$a_{0}$; then player 2, having observed the player 1's move $a_{0}$, makes a
move $b_{0}=y_{0}(a_{0})$. At stage $1$ player 1 makes a move $a_{1}%
=x_{1}(a_{0},b_{0})$ depending on the previous moves $a_{0}$ and $b_{0}$; then
player 2 makes a move $b_{1}=y_{1}(a_{0},b_{0},a_{1})$, and so on. At stage
$t$ ($t\leq N$) player 1 makes a move $a_{t}=x_{t}(a^{t-1},b^{t-1})$ depending
on the sequences of the previous actions
\[
a^{t-1}=(a_{0},a_{1},a_{2},...,a_{t-1})\ \text{and }b^{t-1}=(b_{0},b_{1}%
,b_{2},...,b_{t-1})
\]
up to time $t-1$, and then player 2 makes a move $y_{t}(a^{t},b^{t-1})$
depending on $a^{t}=(a_{0},a_{1},a_{2},...,a_{t})$ and $b^{t-1}=(b_{0}%
,b_{1},b_{2},...,b_{t-1})$. The game terminates at stage $N$, where $N$ is
some given natural number.

A\textit{\ strategy of player 1 }is a sequence\
\[
\xi=\{x_{0},\ x_{1}(a^{0},b^{0}),\ x_{2}(a^{1},b^{1}),\ x_{3}(a^{2}%
,b^{2}),\ ...,\ x_{N}(a^{N-1},b^{N-1})\}
\]
where $x_{0}$ is the initial action of player 1 and $x_{t}(a^{t-1},b^{t-1})$
($1\leq t\leq N$) is a function specifying what move $a_{t}=x_{t}%
(a^{t-1},b^{t-1})$ should be made at stage $t$ given the history
$(a^{t-1},b^{t-1})$ of the previous moves of the players. To specify a
\textit{strategy of player 2} one has to specify a sequence of functions%
\[
\eta=\{y_{0}(a^{0}),\ y_{1}(a^{1},b^{0}),\ y_{2}(a^{2},b^{1}),\ ...,y_{N}%
(a^{N},b^{N-1})\}
\]
indicating what move $b_{t}=y_{t}(a^{t},b^{t-1})$ should be made at stage $t$
given the history $(a^{t},b^{t-1})$ of the previous moves of the players.

The \textit{outcome of the}\emph{\ }\textit{game} $h(\xi,\eta)$ resulting from
the application of the strategies $\xi$ and $\eta$ is described by the whole
history of play%
\[
h(\xi,\eta)=(a^{N},b^{N})=(a_{0},a_{1},a_{2},...,a_{N},b_{0},b_{1}%
,b_{2},...,b_{N}).
\]
\smallskip Once the outcome $h(\xi,\eta)$ of the game is known, the players
get their payoffs $U(\xi,\eta)=u(h(\xi,\eta))$ and $V(\xi,\eta)=v(h(\xi
,\eta))$,\ where $u(h)$ and $v(h)$ are the given \textit{payoff functions}
defined for all histories $h$.

In the course of the proof of Theorem 2, we will establish the existence of
winning and unbeatable strategies having a special structure:\emph{\ }basic
strategies. We call a strategy \textit{basic} if moves of the player using
this strategy depend only on the previous moves of the rival (there is no need
to memorize your own moves). Thus a basic strategy of player 1 is a sequence
$\xi=\{x_{0},x_{1}(b^{0}),x_{2}(b^{1}),...,x_{N}(b^{N-1})\}$ and a basic
strategy of player 2 is a sequence $\eta=$ $\{y_{0}(a^{0}),$ $y_{1}(a^{1}),$
$y_{2}(a^{2}),$ $...,y_{N}(a^{N})\}$.

\textbf{3. }In chess, possible actions/moves of players 1 and 2 (White and
Black) can be identified with positions on the board. When selecting a move,
the player selects a new position. The payoffs, depending on the history of
play, are defined as follows. If White wins then White gets $1$ and Black $0$;
if Black wins then Black gets $1$ and White $0$; in case of a draw, both get
$1/2$. Illegitimate moves (or sequences of moves) lead, by definition, to a
zero payoff for the corresponding player. Winning strategies defined above in
the general context correspond to winning strategies in chess, and unbeatable
ones to those strategies in chess which guarantee at least a draw.

It is assumed that chess is a \textit{finite-stage} game, which is justified
by the following argument. There is a finite number of chess-pieces and a
finite number of squares on the board, hence there is a finite number of
possible positions. The game automatically terminates as a draw if the same
position occurs at least three times, with the same player having to go.
Therefore the game cannot last more than $N$ stages, where $N$ is a
sufficiently large number.

\textbf{4.} For any history of play $h$, define $f(h)=u(h)-v(h)$. Before
proving Theorem 2, we state two auxiliary propositions.

\textbf{Proposition 1.} \textit{A strategy }$\xi$\textit{\ of player 1 is
winning\ (resp. unbeatable) if for any sequence }$b^{N}=(b_{0},...,b_{N})$
\textit{of moves of player 2, }$f(h(\xi,b^{N}))>0$\textit{\ (resp. }%
$f(h(\xi,b^{N}))\geq0$\textit{). A strategy }$\eta$\textit{\ of player 2 is
winning (resp. unbeatable) if for any sequence }$a^{N}=(a_{0},...,a_{N})$
\textit{of moves }of \textit{player 1, }$f(h(a^{N},\eta))<0$\textit{\ (resp.
}$f(h(a^{N},\eta))\leq0$).

\textit{Proof}. The first assertion follows from the fact that for every
strategy profile $(\xi,\eta)$ the outcome $h(\xi,\eta)$ of the game coincides
with $h(\xi,b^{N})$ where $b^{N}=(b_{0},...,b_{N})$ is the sequence of moves
of player 2 generated by the strategy profile $(\xi,\eta)$. The second
assertion is a consequence of the fact that for every strategy profile
$(\xi,\eta)$ the outcome $h(\xi,\eta)$ of the game coincides with
$h(a^{N},\eta)$ where $a^{N}=(a_{0},...,a_{N})$ is the sequence of moves of
player 1 generated by the strategy profile $(\xi,\eta)$.\hfill$\Box$

\textbf{Proposition 2.} \textit{One of the assertions (i) - (iii) listed in
Theorem 2 holds if and only if the following two conditions are satisfied:}

\textit{(I) Either player 1 has a winning strategy, or player 2 has an
unbeatable strategy.}

\textit{(II) Either player 2 has a winning strategy, or player 1 has an
unbeatable strategy.}

\textit{Proof.}\emph{\ }\textit{\textquotedblright If\textquotedblright.}
Suppose (I) and (II) hold. Let us show that one of the assertions (i)-(iii)
holds. If (i) is valid, the assertion is proved. Suppose (i) is not valid.
Then by virtue of (I), player 2 has an unbeatable strategy. If (ii) holds, the
assertion is proved. Suppose not only (i) but also (ii) fails to hold. Then by
virtue of (I) and (II), (iii) holds.

\textit{\textquotedblright Only if\textquotedblright.} If (i) (resp. (ii))
holds, then player 1 (resp. player 2) has a winning, and consequently,
unbeatable strategy, which implies (I) and (II). Similarly, (iii) yields (I)
and (II).\hfill$\Box$

\textbf{5. }\textit{Proof of Theorem 2 }(cf. [2, pp. 147-148] and [3]).
According to the duality principle of first-order logic, the propositions

($\mathbf{P}$) $\exists x_{0}\ \;\forall b_{0}\ \;\exists x_{1}(b^{0}%
)\ \;\forall b_{1}\ \;\exists x_{2}(b^{1})\ \ \;...\;\;\forall b_{N-1}%
\;\;\exists x_{N}(b^{N-1})\ \;\forall b_{N}$ $:$%
\[
f(x_{0},x_{1}(b^{0}),x_{2}(b^{1}),...,x_{N}(b^{N-1}),\;b_{0},b_{1}%
,b_{2},...,b_{N})\ >0
\]
and

($\mathbf{\bar{P}}$) $\forall a_{0}\ \;\exists y_{0}(a^{0})\ \;\forall
a_{1}\ \;\exists y_{1}(a^{1})\ \;\forall a_{2}\ \;\exists y_{2}(a^{2}%
)\ \;...\;\;\forall a_{N}\ \;\exists y_{N}(a^{N})$ $:$%
\[
f(a_{0},a_{1},a_{2},...,a_{N},\;y_{0}(a^{0}),y_{1}(a^{1}),y_{2}(a^{2}%
),...,y_{N}(a^{N}))\leq0
\]
are the \textit{negations} of each other, and so either ($\mathbf{P}$) or
($\mathbf{\bar{P}}$) holds. If ($\mathbf{P}$) holds, then player 1 has a
winning basic strategy $\{x_{0},x_{1}(b^{0}),x_{2}(b^{1}),...,x_{N}%
(b^{N-1})\}$. If ($\mathbf{\bar{P}}$) is valid, then player 2 has an
unbeatable basic strategy $\{y_{0}(a^{0}),y_{1}(a^{1})$, $y_{2}(a^{2})$,...,
$y_{N}(a^{N})$\}. This proves (I). To prove (II) replace in the above argument
$f$ by $-f$.\hfill$\Box$\bigskip

\begin{center}
REFERENCES\bigskip
\end{center}

[1] R. Amir, I. V. Evstigneev and K. R. Schenk-Hopp\'{e}. Asset market games
of survival: A synthesis of evolutionary and dynamic games. \textit{Ann.
Finance }\textbf{9} (2013), 121--144.

[2] A. Kechris. \textit{Classical descriptive set theory}. Springer, 1995.

[3] Yu. Khomskii. \textit{Intensive course on infinite games}. Sofia
University, June 2010.

[4] F. Kojima. Stability and instability of the unbeatable strategy in dynamic
processes. \textit{Int. J. Econ. Theory} \textbf{2} (2006), 41--53.

[5] U. Schwalbe and P. Walker. Zermelo and the early history of game
theory.\ \textit{Games Econ. Behav.} \textbf{34} (2001), 123--137.

[6] E. Zermelo. \"{U}ber eine Anwendung der Mengenlehre auf die Theorie des
Schachspiels. In: E.W. Hobson and A.E.H. Love (eds.), \textit{Proceedings of
the Fifth International Congress of Mathematicians} (Cambridge 1912), Vol. 2.
Cambridge University Press, 1913, pp. 501--504.
\end{document}